\documentclass[11pt]{amsart}
\usepackage{bbm,amsfonts,latexsym,amsmath,amscd,amssymb,fancyhdr}
\usepackage[all]{xy}

\theoremstyle{plain}
\newtheorem{theorem}{Theorem}
\numberwithin{theorem}{section}

\newtheorem{definition}{Definition}
\numberwithin{definition}{section}

\newtheorem{lemma}{Lemma}
\numberwithin{lemma}{section}

\newtheorem{proposition}{Proposition}
\numberwithin{proposition}{section}

\newtheorem{remark}{Remark}
\numberwithin{remark}{section}

\numberwithin{example}{section}

\numberwithin{equation}{section}

\newcommand {\be}{\begin{equation}}
\newcommand {\ee}{\end{equation}}

\newcommand{\h}{\begin{eqnarray*}}
\newcommand{\e}{\end{eqnarray*}}

\newcommand{\CC}{\mathbf{C}}

\newcommand{\ZZ}{\mathbf{Z}}

\newcommand{\ii}{\sqrt{-1}}

\newcommand{\pz}{p_1(TX)}

\newcommand{\px}{p_1(TX)}

\newcommand{\pw}{p_1(W)}

\begin{document}

\title[Modular Forms and Generalized Anomaly Cancellation Formulas]{Modular Forms and  Generalized Anomaly Cancellation Formulas}
\author{Fei Han}
\address{Fei Han, Department of Mathematics, National University of Singapore,
 Block S17, 10 Lower Kent Ridge Road,
Singapore 119076 (mathanf@nus.edu.sg)}
\author{Kefeng Liu}
\address{Kefeng Liu, Department of Mathematics, University of California at Los Angeles,
Los Angeles, CA 90095, USA (liu@math.ucla.edu) and  Center of Mathematical Sciences, Zhejiang University, 310027, P.R. China}
\author{Weiping Zhang}
\address{Weiping Zhang, Chern Institute of Mathematics \& LPMC, Nankai
University, Tianjin 300071, P.R. China. (weiping@nankai.edu.cn)}
\maketitle

\begin{abstract} In this paper, we generalize the anomaly cancellation formulas in \cite{AW, Liu1, HZ2} to  the cases that an auxiliary bundle $W$ as well as a complex line bundle $\xi$ are involved with no conditions on the first Pontryagin forms being assumed. 

\end{abstract}

\section*{Introduction}
In \cite{AW}, gravitational anomaly cancellation formulas are derived from direct computations. In particular, in dimension 12, the Alvarez-Gaum$\mathrm{\acute{e}}$ and Witten ``miraculous cancellation" formula can be written as
\be \{\widehat{L}(TX, \nabla^{TX})\}^{(12)}=\{8\widehat{A}(TX, \nabla^{TX})\mathrm{ch}(T_\CC X, \nabla^{T_\CC X})\}^{(12)}-32\{\widehat{A}(TX, \nabla^{TX})\}^{(12)},\ee where $X$ is a 12 dimensional Riemannian manifold, $\nabla^{TX}$ is the associated Levi-Civita connection, $T_\CC X$ is the complexification of $TX$ (with the induced Hermitian connection $\nabla^{T_\CC X}$) and $\widehat{L}(TX, \nabla^{TX}), \widehat{A}(TX, \nabla^{TX})$ are the Hirzebruch characteristic forms (see (1.1)).

In \cite{Liu1}, Liu generalizes (0.1) to general $8m+4$ dimension by developing modular invariance properties of characteristic forms. Actually, in \cite{Liu1}, Liu obtains a more general cancellation formula by including an auxiliary bundle $W$. More precisely, assume $X$ to be $8m+4$ dimensional and $W$ be a rank $2l$ Euclidean vector bundle over $X$ with a Euclidean connection $\nabla^{W}$ and curvature $R^W=\nabla^{W,2}$, if $p_1(TX, \nabla^{TX})=p_1(W, \nabla^{W})$, then the following identity holds,
\be \begin{split} &\left\{\widehat{A}(TX, \nabla^{TX})\mathrm{det}^{1/2}\left(2\cosh\left(\frac{\ii}{4\pi}R^W\right)\right)
\right\}^{(8m+4)}\\
=&\sum_{r=0}^m2^{l+2m+1-6r}\left\{\widehat{A}(TX, \nabla^{TX})\mathrm{ch}(b_r(T_\CC X, W_\CC, \CC^2)) \right\}^{(8m+4)},\end{split}\ee where $b_r(T_\CC X, W_\CC, \CC^2)$'s are virtual complex vector bundles with connections over $X$ canonically determined by $(TX, \nabla^{TX})$ and $(W, \nabla^W)$. In dimension 12, by direct computation, (0.2) becomes,
\be \begin{split} &\left\{\widehat{A}(TX, \nabla^{TX})\mathrm{det}^{1/2}\left(2\cosh\left(\frac{\ii}{4\pi}R^W\right)\right)\right\}^{(12)}\\
=&2^{l-3}\left\{\widehat{A}(TX, \nabla^{TX})\mathrm{ch}(W_\CC, \nabla^{W_\CC})\right\}^{(12)}-2^{l-2}(l-4)\left\{\widehat{A}(TX, \nabla^{TX})\right\}^{(12)}.\end{split}\ee

When $(TX, \nabla^{TX})=(W, \nabla^{W})$, (0.2) gives,
\be \begin{split} &\frac{1}{8}\left\{\widehat{L}(TX, \nabla^{TX})\right\}^{(8m+4)}\\
=&\sum_{r=0}^m2^{6m-6r}\left\{\widehat{A}(TX, \nabla^{TX})\mathrm{ch}(b_r(T_\CC X, T_\CC X, \CC^2))\right\}^{(8m+4)}.\end{split}\ee
As an application (\cite{Lan}), by the Atiyah-Hirzebruch divisibility \cite{AH}, (0.4) implies the Ochanine divisibility \cite{Och}, which asserts that the signature of an $8k+4$-dimensional smooth closed spin manifold is divisible by 16.

To study higher dimensional Rokhlin congruence, Han and Zhang (\cite{HZ1, HZ2}) extend the ``miraculous cancellation" formulas of Alvarez-Gaum$\mathrm{\acute{e}}$, Witten and Liu to a twisted version where an extra complex line bundle (or equivalently a rank 2 real oriented vector bundle) is involved. More precisely, if $\xi$ is a rank 2 real oriented Euclidean vector bundle equipped with a Euclidean connection $\nabla^{\xi}$ and $c=e(\xi, \nabla^{\xi})$ is the associated Euler form, when $p_1(TX, \nabla^{TX})=p_1(W, \nabla^{W})$, the following identity holds,
\be \begin{split} &\left\{\frac{\widehat{A}(TX, \nabla^{TX})\mathrm{det}^{1/2}\left(2\cosh\left(\frac{\ii}{4\pi}R^W\right)\right)}{\cosh^2\left(\frac{c}{2}\right)}\right\}^{(8m+4)}\\
=&\sum_{r=0}^m2^{l+2m+1-6r}\left\{\widehat{A}(TX, \nabla^{TX})\mathrm{ch}(b_r(T_\CC X, W_\CC, \xi_\CC))\cosh\left(\frac{c}{2}\right)   \right\}^{(8m+4)},\end{split}\ee
where $b_r(T_\CC X, W_\CC, \xi_\CC)$'s are virtual complex vector bundles with connections over $X$ canonically determined by $(TX, \nabla^{TX})$, $(W, \nabla^W)$ and $(\xi, \nabla^{\xi})$. Obviously, when $\xi$ is trivial and $c=0$, (0.5) reduces to (0.2).

When dim$X=12$ and $(TX, \nabla^{TX})=(W, \nabla^{W})$, (0.5) gives,
\be \begin{split} &\left\{\frac{\widehat{L}(TX, \nabla^{TX})}{\cosh^2\left(\frac{c}{2}\right)}\right\}^{(12)}\\
=&\left\{\left[8\widehat{A}(TX, \nabla^{TX})\mathrm{ch}(W_\CC, \nabla^{W_\CC})-32\widehat{A}(TX, \nabla^{TX}) \right.\right.\\
&\left.\left.\ \ \ -24\widehat{A}(TX, \nabla^{TX})(e^c+e^{-c}-2)\right]
\cosh\left(\frac{c}{2}\right)\right\}^{(12)},\end{split}\ee which extends the Alvarez-Gaum$\mathrm{\acute{e}}$ and Witten ``miraculous cancellation" formula (0.1) in dimension 12.

Note that (0.2) and (0.5) only hold under the condition $p_1(TX, \nabla^{TX})=p_1(W, \nabla^{W})$. In this paper, we study what if we remove this condition. We find that the difference between the left hand sides and the right hand sides in (0.2) and (0.5) can actually be written in the form $$(p_1(TX, \nabla^{TX})-p_1(W, \nabla^{W}))\cdot \mathcal{R},$$ where $\mathcal{R}$ is some characteristic forms canonically determined by $(TX, \nabla^{TX})$, $(W, \nabla^W)$ and $(\xi, \nabla^{\xi})$. For example, we find that in dimension 12, the following identity holds (for simplicity, we drop the connections),
\be \begin{split} &\left\{\widehat{A}(TX)\mathrm{det}^{1/2}\left(2\cosh\left(\frac{\ii}{4\pi}R^W\right)\right)\right\}^{(12)}\\
&-2^{l-3}\left\{\widehat{A}(TX)\mathrm{ch}(W_\CC)\right\}^{(12)}+2^{l-2}(l-4)\left\{\widehat{A}(TX)\right\}^{(12)} \\
=&(\px-\pw)\\
&\cdot \left\{\frac{e^{\frac{1}{24}(\px-\pw)}-1}{\px-\pw}\right.\\
&\ \ \ \ \cdot \left.\left[\widehat{A}(TX)\left(2^{l-3}\mathrm{ch}(W_\CC)-2^{l-2}(l-4)\right)
-\widehat{A}(TX)\mathrm{det}^{1/2}\left(2\cosh\left(\frac{\ii}{4\pi}R^W\right)\right)\right]
 \right\}^{(8)}.\end{split}\ee We will give similar general results for $8m+4$ and $8m$ dimensions in Theorem 1.1 and discuss various special cases in Corollaries 1.2-1.5. We obtain our generalized anomaly cancellation formulas still by developing modular invariance of characteristic forms.

To obtain our cancellation formulas, we were also inspired by the Green-Schwarz mechanism.  In \cite{GS}, Green and Schwarz discovered that the anomaly in type I string theory with the gauge group $SO(32)$ cancels because of an extra "classical" contribution from a 2-form field. One key step is that when the gauge group is 496 dimensional, the anomaly $I_{12}$ can be written as (c.f. \cite{S})
\be I_{12}=(p_1(Z)-p_1(F))I_8.\ee Our cancellation formulas in Theorem 1.1 and its corollaries are of same pattern. We hope they can find applications in physics.

\section{Results}
The purpose of this section is to state our main results. We first recall the definitions of some characteristic forms to be used in Section 1.1 and then present our generalized anomaly cancellation formulas in Section 1.2.
\subsection{Some characteristic forms}
Let $X$ be a $4k$ dimensional Riemannian manifold. Let $\nabla^{TX}$ be the associated Levi-Civita connection and
$R^{TX}=\nabla^{TX,\ 2}$ be the curvature
of $\nabla^{TX}$. Let
$\widehat{A}(TX, \nabla^{TX})$ and $\widehat{L}(TX, \nabla^{TX})$ be the
Hirzebruch characteristic forms defined respectively by (cf.
\cite{Z}):
\begin{equation}
\begin{split}
&\widehat{A}(TX, \nabla^{TX}) ={\det}^{1/2}\left({{\sqrt{-1}\over
4\pi}R^{TX} \over \sinh\left({ \sqrt{-1}\over
4\pi}R^{TX}\right)}\right), \\ &\widehat{L}(TX, \nabla^{TX})
={\det}^{1/2}\left({{\sqrt{-1}\over 2\pi}R^{TX} \over \tanh\left({
\sqrt{-1}\over 4\pi}R^{TX}\right)}\right).
\end{split}
\end{equation}

Let $F,G$ be two Hermitian vector bundles over $X$ carrying
Hermitian connections $\nabla^F, \nabla^G$ respectively. Let
$R^F=\nabla^{F,\ 2}$ (resp. $R^G=\nabla^{G,\ 2}$) be the curvature
of $\nabla^F$ (resp. $\nabla^G$). If we set the formal difference
$H=F-G$, then $H$ carries an induced Hermitian connection
$\nabla^H$ in an obvious sense. We define the associated Chern
character form as (cf. \cite{Z})
\be  {\rm ch}(H,\nabla^H)={\rm tr}\left[{\rm exp}\left(\frac{\sqrt{-1}}{2
\pi}R^{F}\right)\right]-{\rm tr}\left[{\rm
exp}\left(\frac{\sqrt{-1}}{2 \pi}R^{G}\right)\right].\ee

For any complex number $t$, let
$$\Lambda_t(F)=\mathbf{C}|_X+tF+t^2\Lambda^2(F)+\cdots ,
\ \  S_t(F)=\mathbf{C}|_X+tF+t^2S^2(F)+\cdots$$  denote respectively
the total exterior and symmetric powers  of $F$, which live in
$K(X)[[t]].$ The following relations between these two operations
(\cite{At}, Chap. 3) hold, \be S_t(F)=\frac{1}{\Lambda_{-t}(F)},\ \ \ \
 \Lambda_t(F-G)=\frac{\Lambda_t(F)}{\Lambda_t(G)}.\ee

The connections
$\nabla^F, \nabla^G$ naturally induce connections on $\Lambda_tF, S_tF
$ etc. Moreover, if $\{\omega_i \}$, $\{{\omega_j}' \}$
are formal Chern roots for Hermitian vector bundles $F$, $G$
respectively, then [\cite{Hir}, Chap. 1] \be
\mathrm{ch}\left(\Lambda_t{(F)},
\nabla^{\Lambda_t(F)}\right)=\prod\limits_i(1+e^{\omega_i}t).\ee

We have the following formulas for Chern character forms,
\be{\rm ch}\left(S_t(F), \nabla^{S_t(F)} \right)=\frac{1}{{\rm
ch}\left(\Lambda_{-t}(F),\nabla^{\Lambda_{-t}(F)}
\right)}=\frac{1}{\prod\limits_i (1-e^{\omega_i}t)}\ ,\ee \be{\rm
ch}\left(\Lambda_t(F-G), \nabla^{\Lambda_t(F-G)} \right)=\frac{{\rm
ch}\left(\Lambda_{t}(F),\nabla^{\Lambda_t(F)} \right)}{{\rm
ch}\left(\Lambda_{t}(G),\nabla^{\Lambda_t(G)}
\right)}=\frac{\prod\limits_i(1+e^{\omega_i}t)}{\prod\limits_j(1+e^{{\omega_j}'}t)}\
.\ee

\subsection{Statement of results}

We make the same assumptions and use the same notations as in Section 1.1.

Let $W$ be a rank $2l$ real Euclidean vector bundle over $X$ carrying a Euclidean connection $\nabla^W$. Let $R^W=\nabla^{W,2}$ be the curvature of $\nabla^W$. If $W$ is spin, let $\Delta(W)=S^+(W)\oplus S^-(W)$ be the spinor bundle of $W$ with the induced connection $\nabla^{\Delta(W)}$. It's not hard to see that $$\mathrm{ch}(\Delta(W), \nabla^{\Delta(W)})=\mathrm{det}^{1/2}\left(2\cosh\left( \frac{\ii}{4\pi}R^W\right)\right).$$

Let $p_1(TX, \nabla^{X})$ and $p_1(W, \nabla^W)$ be the first Pontryagin forms of $(TX, \nabla^{TX})$ and $(W, \nabla^W)$ respectively.

Let $\xi$ be a rank two real oriented Euclidean vector bundle over $X$ carrying a Eucledean connection $\nabla^{\xi}$. Let $c=e(\xi, \nabla^\xi)$ be the Euler form canonically associated to $\nabla^\xi$.

For simplicity, from now on, when there is no ambiguity, we will write characteristic forms without specifying the connections.

In the following, we will define some virtual bundles with connections and some differential forms on $X$ associated to $(TX, \nabla^{TX}), (W, \nabla^W)$ and $(\xi, \nabla^\xi)$.

If $E$ is a vector bundle (real or complex)
over $X$, set $\widetilde{E}=E-{{\rm dim}E}$ in $KO(X)$ or $K(X)$.

If $E$ is a real Euclidean vector bundle over $X$ carrying a
Euclidean connection $\nabla^E$, then its complexification
$E_\CC=E\otimes \CC$ is a complex vector bundle over
$X$ carrying a canonically induced Hermitian metric from that of
$E$, as well as a Hermitian connection $\nabla^{E_\CC}$
induced from $\nabla^E$.

If $\omega$ is a differential form, denote the degree $j$-component of $\omega$ by $\omega^{(j)}$.

Let $q=e^{2\pi \sqrt{-1}\tau}$ with $\tau \in \mathbf{H}$, the upper
half complex plane. Let $T_\CC X$ be the complexification of $TX$.

Set

\be
\begin{split} \Theta_2(T_\CC X, W_\CC, \xi_\CC)=&\bigotimes_{u=1}^\infty
S_{q^u}(\widetilde{T_\CC X}) \otimes \bigotimes_{v=1}^\infty
\Lambda_{-q^{v-{1\over 2}}}(\widetilde{W_\CC}-2\widetilde{\xi_\CC})\\
&\otimes\bigotimes_{r=1}^{\infty}\Lambda_{q^{r-\frac{1}{2}}}(\widetilde{\xi_\CC})
\otimes\bigotimes_{s=1}^{\infty}\Lambda_{q^{s}}(\widetilde{\xi_\CC}) ,\\
\end{split}\ee which is
an element in $K(X)[[q^{1\over2}]]$.

Clearly, $\Theta_2(T_\CC X, W_\CC, \xi_\CC)$ admits
formal Fourier expansion in $q^{1/2}$ as
\be \Theta_2(T_\CC X, W_\CC, \xi_\CC)=B_0(T_\CC X, W_\CC, \xi_\CC)+B_1(T_\CC X, W_\CC, \xi_\CC)q^{1/2}+\cdots,\ee
where the $B_j$'s are elements in the semi-group formally generated by
complex vector bundles over $X$. Moreover, they carry canonically induced
connections denoted by $\nabla^{B_j}$ and let
$\nabla^{\Theta_2}$ be the induced connections
with $q^{1/2}$-coefficients on
$\Theta_2$ from the $\nabla^{B_j}$.

Consider the $q$-series:

\be \delta_1(\tau)=\frac{1}{4}+6\sum_{n=1}^{\infty}\underset{d\ odd}{\underset{d|n}{\sum}}dq^n={1\over 4}+6q+6q^2+\cdots,\ee
\be \varepsilon_1(\tau)=\frac{1}{16}+\sum_{n=1}^{\infty}\underset{d|n}{\sum}(-1)^dd^3q^n={1\over
16}-q+7q^2+\cdots,\ee
\be \delta_2(\tau)=-\frac{1}{8}-3\sum_{n=1}^{\infty}\underset{d\ odd}{\underset{d|n}{\sum}}dq^{n/2}=-{1\over 8}-3q^{1/2}-3q-\cdots,\ee
\be \varepsilon_2(\tau)=\sum_{n=1}^{\infty}\underset{n/d\ odd}{\underset{d|n}{\sum}}d^3q^{n/2}=q^{1/2}+8q+\cdots,\ee and the Eisenstein series
$$E_2(\tau)=1-24\sum_{n=1}^{\infty}\left(\underset{d|n}{\sum}d\right)q^n=1-24q-72q^2-96q^3-\cdots.$$
\begin{remark}$\delta_1$ and $\varepsilon_1$ will only be used later in the proof of our results. We list them here for completeness.
\end{remark}

Now we define the virtual bundles with connections and the differential forms on $X$ associated to $(TX, \nabla^{TX}), (W, \nabla^W)$ and $(\xi, \nabla^\xi)$. This will be done in two cases respectively.

\noindent {\bf Case 1:} $\dim X=8m+4$.

Define virtual complex vector bundles $b_r(T_\CC X, W_\CC, \xi_\CC)$ on $X$, $0\leq r \leq
m,$ via the equality \be \Theta_2(T_\CC X, W_\CC, \xi_\CC)\equiv \sum_{r=0}^m
b_r(8\delta_2)^{2m+1-2r}\varepsilon_2^r\  \ \ \
\mathrm{mod}\  q^{\frac{m+1}{2}}\cdot K(X)[[q^{\frac{1}{2}}]].
\ee

It's not hard to see that each $b_r, 0\leq r \leq m$, is a
canonical integral linear combination of $B_j(T_\CC X, W_\CC, \xi_\CC), 0\leq j\leq r.$
These $b_r$'s carry canonically induced metrics and connections. It's easy to calculate that
\be b_0=-\mathbf{C}, \ b_1=W_\CC-3\xi_\CC+\mathbf{C}^{48m-2l+30}.\ee

Define degree $8m$ differential forms $\beta_r(\nabla^{TX}, \nabla^{W}, \nabla^{\xi})$ on $X$, $0\leq r \leq
m,$ via the equality \be \begin{split}
&\left\{\frac{e^{\frac{1}{24}E_2(\tau)(\pz-\pw)}-1}{\pz-\pw}\widehat{A}(TX)\cosh\left(\frac{c}{2}\right)\mathrm{ch}\left(\Theta_2(T_\CC X, W_\CC, \xi_\CC)\right)\right\}^{(8m)}\\
\equiv &\sum_{r=0}^m
\beta_r(8\delta_2)^{2m+1-2r}\varepsilon_2^r\  \ \ \
\mathrm{mod}\  q^{\frac{m+1}{2}}\cdot \Omega^{8m}(X)[[q^{\frac{1}{2}}]].\end{split}
\ee

It's not hard to see that each $\beta_r, 0\leq r \leq m$, is a
canonical linear combination of degree $8m$ forms of the type $c(\pz-\pw)^a\widehat{A}(TX)\cosh\left(\frac{c}{2}\right)\mathrm{ch}(B_j), 0 \leq j\leq r.$ It's easy to calculate that
\be \begin{split} \beta_0&=-\left\{\frac{e^{\frac{1}{24}(\pz-\pw)}-1}{\pz-\pw}\widehat{A}(TX)\cosh\left(\frac{c}{2}\right)\right\}^{(8m)},\\  \beta_1&=\left\{\frac{e^{\frac{1}{24}(\pz-\pw)}-1}{\pz-\pw}\widehat{A}(TX)\cosh\left(\frac{c}{2}\right)(\mathrm{ch}(W_\CC-3\xi_\CC)
+48m-2l+30)\right\}^{(8m)}.\end{split}\ee

We would like to point out that although $$\beta_0=\left\{\frac{e^{\frac{1}{24}(\pz-\pw)}-1}{\pz-\pw}\widehat{A}(TX)\cosh\left(\frac{c}{2}\right)\mathrm{ch}(b_0)\right\}^{(8m)} $$ and
$$\beta_1=\left\{\frac{e^{\frac{1}{24}(\pz-\pw)}-1}{\pz-\pw}\widehat{A}(TX)\cosh\left(\frac{c}{2}\right)\mathrm{ch}(b_1)\right\}^{(8m)}, $$
generally,
$$\beta_r\neq\left\{\frac{e^{\frac{1}{24}(\pz-\pw)}-1}{\pz-\pw}\widehat{A}(TX)\cosh\left(\frac{c}{2}\right)\mathrm{ch}(b_r)\right\}^{(8m)}, \ r>2.$$

\noindent {\bf Case 2}: $\dim X=8m$.

Define
virtual complex vector bundles $z_r(T_\CC X, W_\CC, \xi_\CC)$ on $X$, $0\leq r \leq
m,$ via the equality \be \Theta_2(T_\CC X, W_\CC, \xi_\CC)\equiv \sum_{r=0}^m
z_r(8\delta_2)^{2m-2r}\varepsilon_2^r\ \ \ \ \
\mathrm{mod}\ q^{\frac{m+1}{2}}\cdot
K(M)[[q^{\frac{1}{2}}]].\ee

Similarly each $z_r(T_\CC X, W_\CC, \xi_\CC), 0\leq r \leq m$, is a
canonical integral linear combination of $B_j(T_\CC X, W_\CC, \xi_\CC), 0\leq j\leq r.$
These $z_r$'s carry canonically induced metrics and connections. It's easy to calculate that
\be z_0=\mathbf{C}, \ z_1=-W_\CC+3\xi_\CC-\mathbf{C}^{48m-2l+6}.\ee

Define degree $8m-4$ differential forms $\zeta_r(\nabla^{TX}, \nabla^{W}, \nabla^{\xi})$ on $X$, $0\leq r \leq
m,$ via the equality \be \begin{split}
&\left\{\frac{e^{\frac{1}{24}E_2(\tau)(\pz-\pw)}-1}{\pz-\pw}\widehat{A}(TX)\cosh\left(\frac{c}{2}\right)\mathrm{ch}\left(\Theta_2(T_\CC X, W_\CC, \xi_\CC)\right)\right\}^{(8m-4)}\\
\equiv &\sum_{r=0}^m
\zeta_r(8\delta_2)^{2m-2r}\varepsilon_2^r\  \ \ \
\mathrm{mod}\  q^{\frac{m+1}{2}}\cdot \Omega^{8m-4}(X)[[q^{\frac{1}{2}}]].\end{split}
\ee

It's not hard to see that each $\zeta_r, 0\leq r \leq m$, is a
canonical linear combination of degree $8m-4$ forms of the type $c(\pz-\pw)^a\widehat{A}(TX)\cosh\left(\frac{c}{2}\right)\mathrm{ch}(B_j), 0 \leq j\leq r.$ It's easy to calculate that
\be \begin{split} \zeta_0&=\left\{\frac{e^{\frac{1}{24}(\pz-\pw)}-1}{\pz-\pw}\widehat{A}(TX)\cosh\left(\frac{c}{2}\right)\right\}^{(8m-4)},\\  \zeta_1&=-\left\{\frac{e^{\frac{1}{24}(\pz-\pw)}-1}{\pz-\pw}\widehat{A}(TX)\cosh\left(\frac{c}{2}\right)(\mathrm{ch}(W_\CC-3\xi_\CC)
+48m-2l+6)\right\}^{(8m-4)}.\end{split}\ee

Similar to $\beta_r$, generally,
$$\zeta_r\neq\left\{\frac{e^{\frac{1}{24}(\pz-\pw)}-1}{\pz-\pw}\widehat{A}(TX)\cosh\left(\frac{c}{2}\right)\mathrm{ch}(z_r)\right\}^{(8m-4)}, \ r>2.$$

We can now state our main theorem as follows.

\begin{theorem} (1)When $\dim X=8m+4$, one has

\be \begin{split} &\left\{\frac{\widehat{A}(TX)\mathrm{det}^{1/2}\left(2\cosh\left(\frac{\ii}{4\pi}R^W\right)\right)}{\cosh^2\left(\frac{c}{2}\right)}\right\}^{(8m+4)}\\
&-\sum_{r=0}^m2^{l+2m+1-6r}\left\{\widehat{A}(TX)\mathrm{ch}(b_r(T_\CC X, W_\CC, \xi_\CC))\cosh\left(\frac{c}{2}\right)   \right\}^{(8m+4)}\\
=&(\px-\pw)\cdot \mathfrak{B}(\nabla^{TX}, \nabla^{W}, \nabla^{\xi}),\end{split}\ee
where
\be \begin{split}
&\mathfrak{B}(\nabla^{TX}, \nabla^{W}, \nabla^{\xi})\\
=&\sum_{r=0}^m2^{l+2m+1-6r}\beta_r(\nabla^{TX}, \nabla^{W}, \nabla^{\xi})\\
&-\left\{\frac{e^{\frac{1}{24}(\px-\pw)}-1}{\px-\pw}\cdot\frac{\widehat{A}(TX)\mathrm{det}^{1/2}\left(2\cosh\left(\frac{\ii}{4\pi}R^W\right)\right)}{\cosh^2\left(\frac{c}{2}\right)}\right\}^{(8m)}.\end{split}\ee

\noindent (2) When $\dim X=8m$, one has
\be \begin{split} &\left\{\frac{\widehat{A}(TX)\mathrm{det}^{1/2}\left(2\cosh\left(\frac{\ii}{4\pi}R^W\right)\right)}{\cosh^2\left(\frac{c}{2}\right)}\right\}^{(8m)}\\
&-\sum_{r=0}^m2^{l+2m-6r}\left\{\widehat{A}(TX)\mathrm{ch}(z_r(T_\CC X, W_\CC, \xi_\CC))\cosh\left(\frac{c}{2}\right)   \right\}^{(8m)}\\
=&(\px-\pw)\cdot \mathfrak{Z}(\nabla^{TX}, \nabla^{W}, \nabla^{\xi}),\end{split} \ee
where
\be \begin{split} &\mathfrak{Z}(\nabla^{TX}, \nabla^{W}, \nabla^{\xi})\\
=&\sum_{r=0}^m2^{l+2m-6r}\zeta_r(\nabla^{TX}, \nabla^{W}, \nabla^{\xi})\\
&-\left\{\frac{e^{\frac{1}{24}(\px-\pw)}-1}{\px-\pw}\cdot\frac{\widehat{A}(TX)\mathrm{det}^{1/2}\left(2\cosh\left(\frac{\ii}{4\pi}R^W\right)\right)}{\cosh^2\left(\frac{c}{2}\right)}\right\}^{(8m-4)}.\end{split}\ee
\end{theorem}
$$ $$

We immediately obtain that

$ \,$\newline
\noindent {\bf Corollary 1.2} (Han-Zhang, \cite{HZ2}) {\it If $p_1(TX, \nabla^{TX})=p_1(W, \nabla^W)$, then

\noindent (1)when $\dim X=8m+4$, the following identity holds,
\be \begin{split} &\left\{\frac{\widehat{A}(TX)\mathrm{det}^{1/2}\left(2\cosh\left(\frac{\ii}{4\pi}R^W\right)\right)}{\cosh^2\left(\frac{c}{2}\right)}\right\}^{(8m+4)}\\
=&\sum_{r=0}^m2^{l+2m+1-6r}\left\{\widehat{A}(TX)\mathrm{ch}(b_r(T_\CC X, W_\CC, \xi_\CC))\cosh\left(\frac{c}{2}\right)   \right\}^{(8m+4)};\end{split}\ee

\noindent (2)when $\dim X=8m$, the following identity holds,
\be \begin{split} &\left\{\frac{\widehat{A}(TX)\mathrm{det}^{1/2}\left(2\cosh\left(\frac{\ii}{4\pi}R^W\right)\right)}{\cosh^2\left(\frac{c}{2}\right)}\right\}^{(8m)}\\
=&\sum_{r=0}^m2^{l+2m-6r}\left\{\widehat{A}(TX)\mathrm{ch}(z_r(T_\CC X, W_\CC, \xi_\CC))\cosh\left(\frac{c}{2}\right)   \right\}^{(8m)}.\end{split}\ee
}

In Corollary 1.2, when $\dim X=8m+4$ and $(W, \nabla^W)=(TX, \nabla^{TX})$ , one has
\be \begin{split} &\frac{1}{8}\left\{\frac{\widehat{L}(TX)}{\cosh^2\left(\frac{c}{2}\right)}\right\}^{(8m+4)}\\
=&\sum_{r=0}^m2^{6m-6r}\left\{\widehat{A}(TX)\mathrm{ch}(b_r(T_\CC X, T_\CC X, \xi_\CC))\cosh\left(\frac{c}{2}\right)   \right\}^{(8m+4)}.\end{split}\ee
This formula is used in \cite{HZ2} to study higher dimensional Rokhlin type congruences.

When $(\xi, \nabla^\xi)=(\mathbf{R}^2, d)$, from Theorem 1.1, we obtain that

$ \,$\newline
\noindent {\bf Corollary 1.3} {\it (1)When $\dim X=8m+4$, one has
\be \begin{split} &\left\{\widehat{A}(TX)\mathrm{det}^{1/2}\left(2\cosh\left(\frac{\ii}{4\pi}R^W\right)\right)
\right\}^{(8m+4)}\\
&-\sum_{r=0}^m2^{l+2m+1-6r}\left\{\widehat{A}(TX)\mathrm{ch}(b_r(T_\CC X, W_\CC, \CC^2)) \right\}^{(8m+4)}\\
=&(\px-\pw)\cdot \mathfrak{B}(\nabla^{TX}, \nabla^{W}, d),\end{split}\ee
where
\be \begin{split}
&\mathfrak{B}(\nabla^{TX}, \nabla^{W}, d)\\
=&\sum_{r=0}^m2^{l+2m+1-6r}\beta_r(\nabla^{TX}, \nabla^{W},d)\\
&-\left\{\frac{e^{\frac{1}{24}(\px-\pw)}-1}{\px-\pw}\cdot\widehat{A}(TX)\mathrm{det}^{1/2}\left(2\cosh\left(\frac{\ii}{4\pi}R^W\right)\right)\right\}^{(8m)}.\end{split}\ee

\noindent (2) When $\dim X=8m$, one has
\be \begin{split} &\left\{\widehat{A}(TX)\mathrm{det}^{1/2}\left(2\cosh\left(\frac{\ii}{4\pi}R^W\right)\right)\right\}^{(8m)}\\
&-\sum_{r=0}^m2^{l+2m-6r}\left\{\widehat{A}(TX)\mathrm{ch}(z_r(T_\CC X, W_\CC, \CC^2))   \right\}^{(8m)}\\
=&(\px-\pw)\cdot \mathfrak{Z}(\nabla^{TX}, \nabla^{W}, d),\end{split} \ee
where
\be \begin{split} &\mathfrak{Z}(\nabla^{TX}, \nabla^{W}, d)\\
=&\sum_{r=0}^m2^{l+2m-6r}\zeta_r(\nabla^{TX}, \nabla^{W}, d)\\
&-\left\{\frac{e^{\frac{1}{24}(\px-\pw)}-1}{\px-\pw}\cdot\widehat{A}(TX)\mathrm{det}^{1/2}\left(2\cosh\left(\frac{\ii}{4\pi}R^W\right)\right)\right\}^{(8m-4)}.\end{split}\ee}

It's interesting to notice that the above anomaly cancellation formulas also imply some integrality results. From Corollary 1.3, we can see that if $X$ is an $8m+4$ dimensional closed spin manifold and $W$ is a $2l$ dimensional spin vector bundle over $X$, then when $l\geq 4m-1$, 
\be \int_X (\px-\pw)\cdot \mathfrak{B}(\nabla^{TX}, \nabla^{W}, d) \ee 
is an integer. Moreover, if $X$ is string, then 
\be \int_X \pw \cdot \mathfrak{B}(\nabla^{TX}, \nabla^{W}, d) \ee  is an integer.  Similarly, we can see that if $X$ is an $8m$ dimensional closed spin manifold and $W$ is a $2l$ dimensional spin vector bundle over $X$, then when $l\geq 4m$, 
\be \int_X (\px-\pw)\cdot \mathfrak{Z}(\nabla^{TX}, \nabla^{W}, d) \ee 
is an integer. Moreover, if $X$ is string, then 
\be \int_X \pw \cdot \mathfrak{Z}(\nabla^{TX}, \nabla^{W}, d) \ee  is an integer.

From Corollary 1.3, we immediately obtain that

$ \,$\newline
\noindent {\bf Corollary 1.4} (Liu, \cite{Liu1}) {\it If $p_1(TX, \nabla^{TX})=p_1(W, \nabla^W)$, then \newline
\noindent (1)When $\dim X=8m+4$, the following identity holds,

\be \begin{split} &\left\{\widehat{A}(TX)\mathrm{det}^{1/2}\left(2\cosh\left(\frac{\ii}{4\pi}R^W\right)\right)
\right\}^{(8m+4)}\\
=&\sum_{r=0}^m2^{l+2m+1-6r}\left\{\widehat{A}(TX)\mathrm{ch}(b_r(T_\CC X, W_\CC, \CC^2)) \right\}^{(8m+4)};\end{split}\ee

\noindent (2) When $\dim X=8m$, one has
\be \begin{split} &\left\{\widehat{A}(TX)\mathrm{det}^{1/2}\left(2\cosh\left(\frac{\ii}{4\pi}R^W\right)\right)\right\}^{(8m)}\\
=&\sum_{r=0}^m2^{l+2m-6r}\left\{\widehat{A}(TX)\mathrm{ch}(z_r(T_\CC X, W_\CC, \CC^2))   \right\}^{(8m)}. \end{split}\ee

}

In Corollary 1.4, when dim$X=8m+4$ and $(W, \nabla^W)=(TX, \nabla^{TX})$ , one has (\cite{Liu1}, \cite{Lan})
\be \begin{split} &\frac{1}{8}\left\{\widehat{L}(TX)\right\}^{(8m+4)}\\
=&\sum_{r=0}^m2^{6m-6r}\left\{\widehat{A}(TX)\mathrm{ch}(b_r(T_\CC X, T_\CC X, \CC^2))\right\}^{(8m+4)}.\end{split}\ee
This formula implies the Ochanine divisibility \cite{Och}, which asserts that the signature of an $8k+4$-dimensional smooth closed spin manifold is divisible by 16.

As examples, we give the explicit formulas when the dimension of $X$ is $4$, $8$ and $12$. Using (1.14), (1.16), (1.18) and (1.20), by direct computations, we have, \newline
$ \,$\newline
\noindent {\bf Corollary 1.5}
{\it (1) when $\dim X=4$, the following identities hold,
\be \begin{split} &\left\{\frac{\widehat{A}(TX)\mathrm{det}^{1/2}\left(2\cosh\left(\frac{\ii}{4\pi}R^W\right)\right)}{\cosh^2\left(\frac{c}{2}\right)}\right\}^{(4)}
+2^{l+1}\left\{\widehat{A}(TX)\cosh\left(\frac{c}{2}\right)\right\}^{(4)}\\
=&-2^{l-3}(\px-\pw),\end{split}\ee
\be \begin{split} &\left\{\widehat{A}(TX)\mathrm{det}^{1/2}\left(2\cosh\left(\frac{\ii}{4\pi}R^W\right)\right)\right\}^{(4)}
+2^{l+1}\left\{\widehat{A}(TX)\right\}^{(4)}\\
=&-2^{l-3}(\px-\pw);\end{split}\ee

\noindent (2) when  $\dim X=8$, the following identities hold,
\be \begin{split} &\left\{\frac{\widehat{A}(TX)\mathrm{det}^{1/2}\left(2\cosh\left(\frac{\ii}{4\pi}R^W\right)\right)}{\cosh^2\left(\frac{c}{2}\right)}\right\}^{(8)}\\
&-\left\{\left[-2^{l-4}\widehat{A}(TX)\mathrm{ch}(W_\CC)+2^{l-3}(l+8)\widehat{A}(TX)+3\cdot 2^{l-4}\widehat{A}(TX)(e^c+e^{-c}-2)\right]
\cosh\left(\frac{c}{2}\right)\right\}^{(8)} \\
=&(\px-\pw)\\
&\cdot \left\{\frac{e^{\frac{1}{24}(\px-\pw)}-1}{\px-\pw}\right.\\
&\ \ \ \ \cdot \left[\widehat{A}(TX)\cosh\left(\frac{c}{2}\right)\left(-2^{l-4}\mathrm{ch}(W_\CC)+2^{l-3}(l+8)+3\cdot 2^{l-4}(e^c+e^{-c}-2)\right)\right.\\
&\ \ \ \ \ \ \ \ \left.\left.-\frac{\widehat{A}(TX)\mathrm{det}^{1/2}\left(2\cosh\left(\frac{\ii}{4\pi}R^W\right)\right)}{\cosh^2\left(\frac{c}{2}\right)}\right] \right\}^{(4)},\end{split}\ee
\be \begin{split} &\left\{\widehat{A}(TX)\mathrm{det}^{1/2}\left(2\cosh\left(\frac{\ii}{4\pi}R^W\right)\right)\right\}^{(8)}\\
&+2^{l-4}\left\{\widehat{A}(TX)\mathrm{ch}(W_\CC)\right\}^{(8)}-2^{l-3}(l+8)\left\{\widehat{A}(TX)\right\}^{(8)} \\
=&(\px-\pw)\\
&\cdot \left\{\frac{e^{\frac{1}{24}(\px-\pw)}-1}{\px-\pw}\right.\\
&\ \ \ \ \cdot \left.\left[\widehat{A}(TX)\left(-2^{l-4}\mathrm{ch}(W_\CC)+2^{l-3}(l+8)\right)
-\widehat{A}(TX)\mathrm{det}^{1/2}\left(2\cosh\left(\frac{\ii}{4\pi}R^W\right)\right)\right]
 \right\}^{(4)};\end{split}\ee

\noindent (3) when $\dim X=12$, the following identities hold,
\be \begin{split} &\left\{\frac{\widehat{A}(TX)\mathrm{det}^{1/2}\left(2\cosh\left(\frac{\ii}{4\pi}R^W\right)\right)}{\cosh^2\left(\frac{c}{2}\right)}\right\}^{(12)}\\
&-\left\{\left[2^{l-3}\widehat{A}(TX)\mathrm{ch}(W_\CC)-2^{l-2}(l-4)\widehat{A}(TX)-3\cdot 2^{l-3}\widehat{A}(TX)(e^c+e^{-c}-2)\right]
\cosh\left(\frac{c}{2}\right)\right\}^{(12)} \\
=&(\px-\pw)\\
&\cdot \left\{\frac{e^{\frac{1}{24}(\px-\pw)}-1}{\px-\pw}\right.\\
&\ \ \ \ \cdot \left[\widehat{A}(TX)\cosh\left(\frac{c}{2}\right)\left(2^{l-3}\mathrm{ch}(W_\CC)-2^{l-2}(l-4)-3\cdot 2^{l-3}(e^c+e^{-c}-2)\right)\right.\\
&\ \ \ \ \ \ \ \ \left.\left.-\frac{\widehat{A}(TX)\mathrm{det}^{1/2}\left(2\cosh\left(\frac{\ii}{4\pi}R^W\right)\right)}{\cosh^2\left(\frac{c}{2}\right)}\right] \right\}^{(8)},\end{split}\ee

\be \begin{split} &\left\{\widehat{A}(TX)\mathrm{det}^{1/2}\left(2\cosh\left(\frac{\ii}{4\pi}R^W\right)\right)\right\}^{(12)}\\
&-2^{l-3}\left\{\widehat{A}(TX)\mathrm{ch}(W_\CC)\right\}^{(12)}+2^{l-2}(l-4)\left\{\widehat{A}(TX)\right\}^{(12)} \\
=&(\px-\pw)\\
&\cdot \left\{\frac{e^{\frac{1}{24}(\px-\pw)}-1}{\px-\pw}\right.\\
&\ \ \ \ \cdot \left.\left[\widehat{A}(TX)\left(2^{l-3}\mathrm{ch}(W_\CC)-2^{l-2}(l-4)\right)
-\widehat{A}(TX)\mathrm{det}^{1/2}\left(2\cosh\left(\frac{\ii}{4\pi}R^W\right)\right)\right]
 \right\}^{(8)}.\end{split}\ee}
 
 \begin{remark} It's not hard to see that (1.41)-(1.44) are respectively equivalent to the following identites,
\be \begin{split} &\left\{e^{\frac{1}{24}(\px-\pw)}\right.\\
&\ \ \ \ \cdot \left[\frac{\widehat{A}(TX)\mathrm{det}^{1/2}\left(2\cosh\left(\frac{\ii}{4\pi}R^W\right)\right)}{\cosh^2\left(\frac{c}{2}\right)}\right.\\
&\ \ \ \ \ \ \ \ \left.\left.-\widehat{A}(TX)\cosh\left(\frac{c}{2}\right)\left(-2^{l-4}\mathrm{ch}(W_\CC)+2^{l-3}(l+8)+3\cdot 2^{l-4}(e^c+e^{-c}-2)\right)\right] \right\}^{(8)}=0,\end{split}\ee

\be \begin{split} 
&\left\{e^{\frac{1}{24}(\px-\pw)}\right.\\
&\ \ \ \ \cdot \left.\left[\widehat{A}(TX)\mathrm{det}^{1/2}\left(2\cosh\left(\frac{\ii}{4\pi}R^W\right)\right)-\widehat{A}(TX)\left(-2^{l-4}\mathrm{ch}(W_\CC)+2^{l-3}(l+8)\right)
\right]
 \right\}^{(8)}=0;\end{split}\ee  

\be \begin{split} 
&\left\{e^{\frac{1}{24}(\px-\pw)}\right.\\
&\ \ \ \ \cdot \left[-\frac{\widehat{A}(TX)\mathrm{det}^{1/2}\left(2\cosh\left(\frac{\ii}{4\pi}R^W\right)\right)}{\cosh^2\left(\frac{c}{2}\right)}\right.\\
&\ \ \ \ \ \ \ \ \left.\left.-\widehat{A}(TX)\cosh\left(\frac{c}{2}\right)\left(2^{l-3}\mathrm{ch}(W_\CC)-2^{l-2}(l-4)-3\cdot 2^{l-3}(e^c+e^{-c}-2)\right)\right] \right\}^{(12)}=0,\end{split}\ee

\be \begin{split} 
&\left\{e^{\frac{1}{24}(\px-\pw)}\right.\\
&\ \ \ \ \cdot \left.\left[\widehat{A}(TX)\mathrm{det}^{1/2}\left(2\cosh\left(\frac{\ii}{4\pi}R^W\right)\right)-\widehat{A}(TX)\left(2^{l-3}\mathrm{ch}(W_\CC)-2^{l-2}(l-4)\right)
\right]
 \right\}^{(12)}=0.\end{split}\ee

These formulas are simply in a form of the multiplication of  $e^{\frac{1}{24}(\px-\pw)}$ to the original anomaly cancellation formulas in dimension 8 and 12 holding under the condition $p_1(X)=p_1(W)$. 

However as pointed out on page 7 and page 8 about the patterns of $\beta_r$ and $\zeta_r$ for $r\geq 2$, we know that the higher ($>12$) dimensional anomaly cancellation formulas are not as simple as the above lower anomaly cancellation formulas, i.e. they are not  simply in a form of the multiplication of  $e^{\frac{1}{24}(\px-\pw)}$ to the original anomaly cancellation formulas holding under the condition $p_1(X)=p_1(W)$.  

 \end{remark}

\section{Proofs}
In this section, we give the proof of Theorem 1.1. To prepare for the proof in Section 2.2, we will first recall some basic knowledge about the Jacobi theta functions, modular forms and Eisenstein series in Section 2.1.

\subsection{Preliminaries} Let $$ SL_2(\mathbf{Z}):= \left\{\left.\left(\begin{array}{cc}
                                      a&b\\
                                      c&d
                                     \end{array}\right)\right|a,b,c,d\in\mathbf{Z},\ ad-bc=1
                                     \right\}
                                     $$
 as usual be the modular group. Let
$$S=\left(\begin{array}{cc}
      0&-1\\
      1&0
\end{array}\right), \ \ \  T=\left(\begin{array}{cc}
      1&1\\
      0&1
\end{array}\right)$$
be the two generators of $ SL_2(\mathbf{Z})$. Their actions on
$\mathbf{H}$ are given by
$$ S:\tau\rightarrow-\frac{1}{\tau}, \ \ \ T:\tau\rightarrow\tau+1.$$

The four Jacobi theta functions are defined as follows (cf.
\cite{C}): \h \theta(v,\tau)=2q^{1/8}\sin(\pi v)
\prod_{j=1}^\infty\left[(1-q^j)(1-e^{2\pi \sqrt{-1}v}q^j)(1-e^{-2\pi
\sqrt{-1}v}q^j)\right]\ ,\e \h \theta_1(v,\tau)=2q^{1/8}\cos(\pi
v)
 \prod_{j=1}^\infty\left[(1-q^j)(1+e^{2\pi \sqrt{-1}v}q^j)
 (1+e^{-2\pi \sqrt{-1}v}q^j)\right]\ ,\e
\h \theta_2(v,\tau)=\prod_{j=1}^\infty\left[(1-q^j)
 (1-e^{2\pi \sqrt{-1}v}q^{j-1/2})(1-e^{-2\pi \sqrt{-1}v}q^{j-1/2})\right]\
 ,\e
\h \theta_3(v,\tau)=\prod_{j=1}^\infty\left[(1-q^j) (1+e^{2\pi
\sqrt{-1}v}q^{j-1/2})(1+e^{-2\pi \sqrt{-1}v}q^{j-1/2})\right]\ .\e
They are all holomorphic functions for $(v,\tau)\in \mathbf{C \times
H}$, where $\mathbf{C}$ is the complex plane and $\mathbf{H}$ is the
upper half plane.

If we act theta-functions by $S$ and $T$, the theta functions obey
the following transformation laws (cf. \cite{C}), \be
\theta(v,\tau+1)=e^{\pi \sqrt{-1}\over 4}\theta(v,\tau),\ \ \
\theta\left(v,-{1}/{\tau}\right)={1\over\sqrt{-1}}\left({\tau\over
\sqrt{-1}}\right)^{1/2} e^{\pi\sqrt{-1}\tau v^2}\theta\left(\tau
v,\tau\right)\ ;\ee \be \theta_1(v,\tau+1)=e^{\pi \sqrt{-1}\over
4}\theta_1(v,\tau),\ \ \
\theta_1\left(v,-{1}/{\tau}\right)=\left({\tau\over
\sqrt{-1}}\right)^{1/2} e^{\pi\sqrt{-1}\tau v^2}\theta_2(\tau
v,\tau)\ ;\ee \be\theta_2(v,\tau+1)=\theta_3(v,\tau),\ \ \
\theta_2\left(v,-{1}/{\tau}\right)=\left({\tau\over
\sqrt{-1}}\right)^{1/2} e^{\pi\sqrt{-1}\tau v^2}\theta_1(\tau
v,\tau)\ ;\ee \be\theta_3(v,\tau+1)=\theta_2(v,\tau),\ \ \
\theta_3\left(v,-{1}/{\tau}\right)=\left({\tau\over
\sqrt{-1}}\right)^{1/2} e^{\pi\sqrt{-1}\tau v^2}\theta_3(\tau
v,\tau)\ .\ee

\begin{definition} Let $\Gamma$ be a subgroup of $SL_2(\mathbf{Z}).$ A modular form over $\Gamma$ is a holomorphic function $f(\tau)$ on $\mathbf{H}\cup
\{\infty\}$ such that for any
 $$ g=\left(\begin{array}{cc}
             a&b\\
             c&d
             \end{array}\right)\in\Gamma\ ,$$
 the following property holds
 $$f(g\tau):=f\left(\frac{a\tau+b}{c\tau+d}\right)=\chi(g)(c\tau+d)^lf(\tau), $$
 where $\chi:\Gamma\rightarrow\mathbf{C}^*$ is a character of
 $\Gamma$ and $l$ is called the weight of $f$.
 \end{definition}

Let
\be E_{2k}=1-\frac{4k}{B_{2k}}\sum_{n=1}^{\infty}\left(\underset{d|n}{\sum}d^{2k-1}\right)q^n \ee
be the Eisenstein series, where $B_{2k}$ is the $2k$-th Bernoulli number.

When $k>1$, $E_{2k}$ is a modular form of weight $2k$ over $SL_2(\mathbf{Z})$. However, unlike other Eisenstein theories, $E_2(\tau)$ is not a modular form over $SL(2,\ZZ)$, instead $E_2(\tau)$ is a quasimodular form over $SL(2,\ZZ)$, satisfying:
\be E_2\left(\frac{a\tau+b}{c\tau+d}\right)=(c\tau+d)^2E_2(\tau)-\frac{6\ii c(c\tau+d)}{\pi}. \ee
In particular, we have
\be E_2(\tau+1)=E_2(\tau),\ee
\be E_2\left(-\frac{1}{\tau}\right)=\tau^2E_2(\tau)-\frac{6\ii\tau}{\pi}.\ee
For the precise definition of quasimodular forms, see \cite{KZ}. 

In the following, let's review some level 2 modular forms.

Let
$$ \Gamma_0(2)=\left\{\left.\left(\begin{array}{cc}
a&b\\
c&d
\end{array}\right)\in SL_2(\mathbf{Z})\right|c\equiv0\ \ (\rm mod \ \ 2)\right\},$$

$$ \Gamma^0(2)=\left\{\left.\left(\begin{array}{cc}
a&b\\
c&d
\end{array}\right)\in SL_2(\mathbf{Z})\right|b\equiv0\ \ (\rm mod \ \ 2)\right\}$$
be the two modular subgroups of $SL_2(\mathbf{Z})$. It is known
that the generators of $\Gamma_0(2)$ are $T,ST^2ST$ and the generators
of $\Gamma^0(2)$ are $STS,T^2STS$.(cf. \cite{C}).

Writing simply
$\theta_j=\theta_j(0,\tau),\ 1\leq j \leq 3,$ we have (cf. \cite{HBJ} and \cite{Liu3}),
$$ \delta_1(\tau)=\frac{1}{8}(\theta_2^4+\theta_3^4), \ \ \ \
\varepsilon_1(\tau)=\frac{1}{16}\theta_2^4 \theta_3^4\ ,$$
$$\delta_2(\tau)=-\frac{1}{8}(\theta_1^4+\theta_3^4), \ \ \ \
\varepsilon_2(\tau)=\frac{1}{16}\theta_1^4 \theta_3^4\ .$$

If $\Gamma$ is a modular subgroup, let
$M_\mathbf{R}(\Gamma)$ denote the ring of modular forms
over $\Gamma$ with real Fourier coefficients.
\begin{lemma} [\protect cf. \cite{Liu1}] One has that $\delta_1(\tau)\ (resp.\ \varepsilon_1(\tau) ) $
is a modular form of weight $2 \ (resp.\ 4)$ over $\Gamma_0(2)$,
$\delta_2(\tau) \ (resp.\ \varepsilon_2(\tau))$ is a modular form
of weight $2\ (resp.\ 4)$ over $\Gamma^0(2)$ and moreover
$M_\mathbf{R}(\Gamma^0(2))=\mathbf{R}[\delta_2(\tau),
\varepsilon_2(\tau)]$. Moreover, we have
transformation laws \be
\delta_2\left(-\frac{1}{\tau}\right)=\tau^2\delta_1(\tau),\ \ \ \ \
\ \ \ \ \
\varepsilon_2\left(-\frac{1}{\tau}\right)=\tau^4\varepsilon_1(\tau).\ee

\end{lemma}

\subsection{Proof of Theorem 1.1}

Set
\be \begin{split} \Theta_1(T_\CC X, W_\CC, \xi_\CC)=&\bigotimes_{u=1}^\infty S_{q^u}(\widetilde{T_\CC X})
\otimes \bigotimes_{v=1}^\infty \Lambda_{q^v}(\widetilde{W_\CC}-2\widetilde{\xi}_\CC)\\
&\otimes \bigotimes_{r=1}^\infty \Lambda_{q^{r-1/2}}(\widetilde{\xi}_\CC)\otimes \bigotimes_{s=1}^\infty \Lambda_{-q^{s-1/2}}(\widetilde{\xi}_\CC).\end{split} \ee

$\Theta_1(T_\CC X, W_\CC, \xi_\CC)$ admits
formal Fourier expansion in $q^{1/2}$ as
\be \Theta_1(T_\CC X, W_\CC, \xi_\CC)=A_0(T_\CC X, W_\CC, \xi_\CC)+ A_1(T_\CC X, W_\CC, \xi_\CC)q^{1/2}+\cdots, \ee
where the $A_j$'s are elements in the semi-group formally generated by
complex vector bundles over $X$. Moreover, they carry canonically induced
connections denoted by $\nabla^{A_j}$, and let
$\nabla^{\Theta_1}$ be the induced connections
with $q^{1/2}$-coefficients on
$\Theta_1$ from the $\nabla^{A_j}$.

To prove part 1 of Theorem 1.1 ($8m+4$-dimensional case), set
\be \begin{split} P_1(\tau):=&\left\{e^{\frac{1}{24}E_2(\tau)(p_1(TX)-p_1(W))}\right.\\
&\ \ \left. \cdot \frac{\widehat{A}(TX)\mathrm{det}^{1/2}\left(2\cosh\left(\frac{\ii}{4\pi}R^W\right)\right)}{\cosh^2\left(\frac{c}{2}\right)}
\mathrm{ch}\left(\Theta_1(T_{\mathbf C}X,W_{\mathbf C}, \xi_\CC)\right)\right\}^{(8m+4)},\end{split}\ee

\be P_2(\tau):=\left\{\widehat{A}(TX)\cosh\left(\frac{c}{2}\right)\mathrm{ch}\left(\Theta_2(T_{\mathbf C}X,W_{\mathbf C}, \xi_\CC)\right)\right\}^{(8m+4)}\ee and

\be \begin{split} \Xi_2(\tau):=
&\left\{\frac{e^{\frac{1}{24}E_2(\tau)(p_1(TX)-p_1(W))}-1}{p_1(TX)-p_1(W)}\right.\\
&\ \ \cdot \left.\widehat{A}(TX)\cosh\left(\frac{c}{2}\right)\mathrm{ch}\left(\Theta_2(T_\CC X, W_\CC, \xi_\CC)\right)\right\}^{(8m)}.\end{split}\ee

We have
\begin{proposition} $P_1(\tau)$ is a modular form of weight $4m+2$ over $\Gamma_0(2)$ while $P_2(\tau)+(p_1(TX)-p_1(W))\Xi_2(\tau)$ is a modular form of weight $4m+2$ over $\Gamma^0(2)$. Moreover, the following identity holds,
\be P_1\left(-\frac{1}{\tau}\right)=2^l\tau^{4m+2}(P_2(\tau)+(p_1(TX)-p_1(W))\Xi_2(\tau)).\ee

\end{proposition}

\begin{proof} Let $\{\pm 2\pi \ii y_k\}$ (resp. $\{\pm 2\pi \ii x_j\}$) be the formal Chern roots for $(W_\CC, \nabla^{W_\CC})$ (resp. $(TM_\CC, \nabla^{TM_\CC})$). Let $c=2\pi \ii u$.

By the Chern root algorithm, we have
\be \begin{split}
P_1(\tau)=2^{l}&\left\{e^{\frac{1}{24}E_2(\tau)(p_1(TX)-p_1(W))}\left(\prod_{j=1}^{4m+2}\left(x_j\frac{\theta'(0,\tau)}{\theta(x_j,\tau)}\right)
\right)\left(\prod_{k=1}^{l}\frac{\theta_{1}(y_k,\tau)}{\theta_{1}(0,\tau)}\right)\right.\\
&\left.\ \ \cdot \frac{\theta_1^2(0,\tau)}{\theta_1^2(u,\tau)}\frac{\theta_3(u,\tau)}{\theta_3(0,\tau)}\frac{\theta_2(u,\tau)}{\theta_2(0,\tau)}\right\}^{(8m+4)},
\end{split}\ee
and
\be \begin{split} &P_2(\tau)+(p_1(TX)-p_1(W))\Xi_2(\tau)\\
=& \left\{e^{\frac{1}{24}E_2(\tau)(p_1(TX)-p_1(W))}\widehat{A}(TX)\cosh\left(\frac{c}{2}\right)\mathrm{ch}\left(\Theta_2(T_\CC X, W_\CC, \xi_\CC)\right)\right\}^{(8m+4)}\\
=&\left\{e^{\frac{1}{24}E_2(\tau)(p_1(TX)-p_1(W))}\left(\prod_{j=1}^{4m+2}\left(x_j\frac{\theta'(0,\tau)}{\theta(x_j,\tau)}\right)
\right)\left(\prod_{j=1}^{l}\frac{\theta_{2}(y_j,\tau)}{\theta_{2}(0,\tau)}\right)\right.\\
&\ \ \left.\cdot\frac{\theta_2^2(0,\tau)}{\theta_2^2(u,\tau)}\frac{\theta_3(u,\tau)}{\theta_3(0,\tau)}\frac{\theta_1(u,\tau)}{\theta_1(0,\tau)}\right\}^{(8m+4)}.\end{split}\ee

Then we can apply the transformation laws (2.1)-(2.4) for theta functions as well as the transformation laws (2.7)-(2.8) to (2.16) and (2.17) to get the desired results.
\end{proof}

We can now proceed to prove part 1 of Theorem 1.1 as follows.

Combining Lemma 2.1 and Proposition 2.1, we can write
\be \begin{split} &P_2(\tau)+(p_1(TX)-p_1(W))\Xi_2(\tau)\\
=&h_0(8\delta_2)^{2m+1}+h_1(8\delta_2)^{2m-1}\varepsilon_2+\cdots+h_m(8\delta_2)\varepsilon_2^m,
\end{split}
\ee
where $h_r\in \Omega^{8m+4}(X), 0\leq r\leq m.$

By the definitions of $b_r(T_\CC X, W_\CC, \xi_\CC)$ and $\beta_r(\nabla^{TX}, \nabla^{W}, \nabla^{\xi})$, it's easy to see that for $0\leq r\leq m$,
\be h_r=\{\widehat{A}(TX)\cosh\left(\frac{c}{2}\right)\mathrm{ch}(b_r(T_\CC X, W_\CC, \xi_\CC))\}^{(8m+4)}+(p_1(TX)-p_1(W))\beta_r(\nabla^{TX}, \nabla^{W}, \nabla^{\xi}).\ee

Therefore (simply denote $b_r(T_\CC X, W_\CC, \xi_\CC)$ and $\beta_r(\nabla^{TX}, \nabla^{W}, \nabla^{\xi})$ by $b_r$ and $\beta_r$),
\be \begin{split} &P_2(\tau)+(p_1(TX)-p_1(W))\Xi_2(\tau)\\
=&\sum_{r=0}^{m}\left(\left\{\widehat{A}(TX)\cosh\left(\frac{c}{2}\right)\mathrm{ch}(b_r)\right\}^{(8m+4)}+(p_1(TX)-p_1(W))\beta_r\right)(8\delta_2)^{2m+1-r}\varepsilon_2^r.
\end{split}
\ee

By (2.9) and (2.15), we have
\be \begin{split}
&P_1(\tau)\\
=&\frac{2^l}{\tau^{4m+2}}\left[P_2\left(-\frac{1}{\tau}\right)+(p_1(TX)-p_1(W))\Xi_2\left(-\frac{1}{\tau}\right)\right]\\
=&\frac{2^l}{\tau^{4m+2}}\left[\left(\left\{\widehat{A}(TX)\cosh\left(\frac{c}{2}\right)\mathrm{ch}(b_0)\right\}^{(8m+4)}+(p_1(TX)-p_1(W))\beta_0\right)\left(8\delta_2\left(-\frac{1}{\tau}\right)\right)^{2m+1}\right.\\
&\ \ \ \ \ \ \ \ \ \ +\cdots\\
&\ \ \ \ \ \ \ \ \ \ +\left(\left\{\widehat{A}(TX)\cosh\left(\frac{c}{2}\right)\mathrm{ch}(b_m)\right\}^{(8m+4)}+(p_1(TX)-p_1(W))\beta_m\right)\\
&\left.\ \ \ \ \ \ \ \ \ \ \ \ \ \ \ \cdot 8\delta_2\left(-\frac{1}{\tau}\right)\left(\varepsilon_2\left(-\frac{1}{\tau}\right)\right)^m\right]\\
=&2^l\left[\left(\left\{\widehat{A}(TX)\cosh\left(\frac{c}{2}\right)\mathrm{ch}(b_0)\right\}^{(8m+4)}+(p_1(TX)-p_1(W))\beta_0\right)\left(8\delta_1\right)^{2m+1}\right.\\
&\ \ \ \ \ +\cdots\\
&\ \ \ \ \ +\left.\left(\left\{\widehat{A}(TX)\cosh\left(\frac{c}{2}\right)\mathrm{ch}(b_m)\right\}^{(8m+4)}+(p_1(TX)-p_1(W))\beta_m\right) (8\delta_1)\varepsilon_1^m\right].\end{split}\ee

Comparing the constant terms of both sides of (2.21), one has
\be \begin{split} &e^{\frac{1}{24}(p_1(TX)-p_1(W))}\frac{\widehat{A}(TX)\mathrm{det}^{1/2}\left(2\cosh\left(\frac{\ii}{4\pi}R^W\right)\right)}{\cosh^2\left(\frac{c}{2}\right)}\\
=&\sum_{r=0}^m 2^{l+2m+1-6r}\left(\left\{\widehat{A}(TX)\cosh\left(\frac{c}{2}\right)\mathrm{ch}(b_r)\right\}^{(8m+4)}+(p_1(TX)-p_1(W))\beta_r\right).  \end{split}\ee So we have
\be \begin{split} &\left\{\frac{\widehat{A}(TX)\mathrm{det}^{1/2}\left(2\cosh\left(\frac{\ii}{4\pi}R^W\right)\right)}{\cosh^2\left(\frac{c}{2}\right)}\right\}^{(8m+4)}\\
&-\sum_{r=0}^m2^{l+2m+1-6r}\left\{\widehat{A}(TX)\mathrm{ch}(b_r(T_\CC X, W_\CC, \xi_\CC))\cosh\left(\frac{c}{2}\right)   \right\}^{(8m+4)}\\
=&(\px-\pw)\cdot \mathfrak{B}(\nabla^{TX}, \nabla^{W}, \nabla^{\xi}),\end{split}\ee
where
\begin{equation*} \begin{split}
&\mathfrak{B}(\nabla^{TX}, \nabla^{W}, \nabla^{\xi})\\
=&\sum_{r=0}^m2^{l+2m+1-6r}\beta_r(\nabla^{TX}, \nabla^{W}, \nabla^{\xi})\\
&-\left\{\frac{e^{\frac{1}{24}(\px-\pw)}-1}{\px-\pw}\cdot\frac{\widehat{A}(TX)\mathrm{det}^{1/2}\left(2\cosh\left(\frac{\ii}{4\pi}R^W\right)\right)}{\cosh^2\left(\frac{c}{2}\right)}\right\}^{(8m)}.\end{split}\end{equation*}

To prove part 2 of Theorem 1.1  ($8m$-dimensional case), set
\be \begin{split} Q_1(\tau):=&\left\{e^{\frac{1}{24}E_2(\tau)(p_1(TX)-p_1(W))}\right.\\
&\ \ \left. \cdot \frac{\widehat{A}(TX)\mathrm{det}^{1/2}\left(2\cosh\left(\frac{\ii}{4\pi}R^W\right)\right)}{\cosh^2\left(\frac{c}{2}\right)}
\mathrm{ch}\left(\Theta_1(T_{\mathbf C}X,W_{\mathbf C}, \xi_\CC)\right)\right\}^{(8m)},\end{split}\ee

\be Q_2(\tau):=\left\{\widehat{A}(TX)\cosh\left(\frac{c}{2}\right)\mathrm{ch}\left(\Theta_2(T_{\mathbf C}X,W_{\mathbf C}, \xi_\CC)\right)\right\}^{(8m)}\ee and

\be \begin{split} \Pi_2(\tau):=
&\left\{\frac{e^{\frac{1}{24}E_2(\tau)(p_1(TX)-p_1(W))}-1}{p_1(TX)-p_1(W)}\right.\\
&\ \ \cdot \left.\widehat{A}(TX)\cosh\left(\frac{c}{2}\right)\mathrm{ch}\left(\Theta_2(T_\CC X, W_\CC, \xi_\CC)\right)\right\}^{(8m-4)}.\end{split}\ee

Similar to the $8m+4$ dimensional case, one has
\begin{proposition}$Q_1(\tau)$ is a modular form of weight $4m$ over $\Gamma_0(2)$ while $Q_2(\tau)+(p_1(TX)-p_1(W))\Pi_2(\tau)$ is a modular form of weight $4m$ over $\Gamma^0(2)$. Moreover, the following identity holds,
\be Q_1\left(-\frac{1}{\tau}\right)=2^l\tau^{4m}(Q_2(\tau)+(p_1(TX)-p_1(W))\Pi_2(\tau)).\ee
\end{proposition}

Then one can prove part 2 of Theorem 1.1 by adopting similar idea as in the above proof of part 1 of Theorem 1.1.

$$ $$
\noindent {\bf Acknowledgements} The first author is partially supported by a start-up grant from National University of Singapore. The second author is partially supported by NSF. The third author is partially supported by  MOE and NNSFC.

\end{document}